\documentclass[a4paper]{amsart}
\usepackage{amsmath}
\usepackage{amsfonts}
\usepackage{amssymb}
\usepackage{latexsym}
\usepackage{epsfig}

\def\QED{~\hfill\nobreak {\textbf{q.e.d.}}}

\newcommand{\lqq}{\lq\lq}
\newcommand{\rqq}{\rq\rq}

\newcommand{\begriff}[1]{\textbf{#1}}

\newtheorem{theorem}{Theorem}
\newtheorem{proposition}{Proposition}
\newtheorem{lemma}{Lemma}
\newtheorem{corollary}{Corollary}

\theoremstyle{definition}

\theoremstyle{remark}

\begin{document}
\title{Regular Flip Equivalence of Surface Triangulations}
\author{Simon A. King}
\address{Simon A. King\\
Department of Mathematics\\
Darmstadt University of Technology\\
Schlossgartenstr.~7\\
64289 Darmstadt\\
  Germany}
\email{king@mathematik.tu-darmstadt.de}
%
%
%

\begin{abstract}
  \noindent
  Any two triangulations of a closed
  surface with the same number of vertices can be transformed into
  each other  by a sequence of flips, provided the
  number of vertices exceeds a number $N$ depending on the
  surface. Examples show that in general $N$ is bigger than the
  minimal number of vertices of a triangulation. The existence
  of $N$ was known, but no estimate. This paper provides an estimate
  for $N$ that is linear in the Euler characteristic of the surface. 
\end{abstract}

\maketitle

\section{Results on flip equivalence}
\label{sec:def}

Let $F$ be a closed surface and let $\chi(F)$ be its Euler
characteristic. A \begriff{singular triangulation} of $F$ is a
graph $T$ embedded in $F$ such that each face of $F\setminus T$
is bounded by an edge path of length three. 
We denote by $v(T)$, $e(T)$ and
$f(T)$ the number of vertices, edges and faces of $T$.  
If $T$ is without loops and multiple edges and has  more than three
faces, then $T$ corresponds to a triangulation of $F$ in the classical
meaning of the word; in order to avoid confusions, we use for it the term
\begriff{regular  triangulation}. 

\begin{figure}[t]
  \begin{center}
    \leavevmode
    \epsfig{file=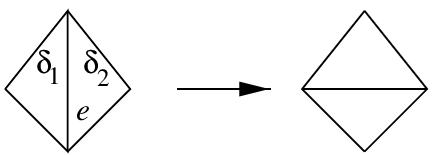}
    \caption{A flip}
    \label{fig:flip}
  \end{center}
\end{figure}

Let $e$ be an edge of a singular triangulation $T$ and suppose that there
are two distinct faces $\delta_1$ and
$\delta_2$ adjacent to $e$. The faces $\delta_1$ and $\delta_2$ 
form a (possibly degenerate) quadrilateral, containing $e$ as a
diagonal. A \begriff{flip} of 
$T$ along $e$ replaces $e$ by the opposite diagonal of this
quadrilateral, see Figure~\ref{fig:flip}.
The flip is called \begriff{regular}, if both $T$ and the result of
the flip are regular triangulations. Two singular (resp.\ regular)
triangulations $T_1$, $T_2$ of a closed surface are called flip
equivalent (resp.\ regularly flip equivalent), if they are related by
a finite sequence of flips (resp.\ regular flips) and isotopy.

The following result is well known, and there are many proofs for
it. There are interesting
applications to the automatic structure of mapping class groups,
see~\cite{Mosher95} or~\cite{RourkeWiest}. 

\begin{lemma}\label{prop:flips}
   Any two singular triangulations $T_1$ and $T_2$ of a closed surface
   $F$ with $v(T_1)=v(T_2)$ are flip equivalent.\qed
\end{lemma}

One might ask whether any two {regular}
triangulations of $F$ with the same number of vertices are
\emph{regularly} flip equivalent. The answer is
\lqq Yes\rqq\ in special cases:
any two regular triangulations of the sphere~\cite{Wagner36}, the
torus~\cite{Dewdney73}, the projective plane 
or the Klein bottle~\cite{NegamiWatanabe90} with the same number of
vertices are regularly flip equivalent. But in general, the answer is
\lqq No\rqq:  it is known that there are 59 different
triangulations of the closed oriented surface of genus six based on the 
complete graph with 12 vertices, see~\cite{Altshuler96}. Such a
triangulation does not admit any regular flip, thus the different
triangulations are not regularly flip equivalent.

This paper is devoted to the proof of the following theorem. A
preliminary version of this paper has been appeared in~\cite{king}. 

\begin{theorem}\label{thm:main}
  Let $F$ be a closed surface  and
  $N(F)= 9450 - 6020\chi(F)$. Any two regular triangulations $T_1$ and
  $T_2$ of $F$ with $v(T_1)= v(T_2)\ge N(F)$ are regularly flip
  equivalent. 
\end{theorem}

Negami~\cite{Negami94} stated the mere existence of $N(F)$
without an estimate.  The estimate in Theorem~\ref{thm:main} 
is far from being best possible, at least for the surfaces up to genus
one. The number $N(F)$ is
negative if and only if $F$ is a sphere, in which case the statement
is true since the transformation by regular flips is always possible,
by Wagner's Theorem~\cite{Wagner36}. We assume in the 
following that $F$ is not the sphere.

\section{Proof of the Theorem}
\label{sec:proof}

\begin{figure}[htb]
  \begin{center}
    \leavevmode
    \epsfig{file=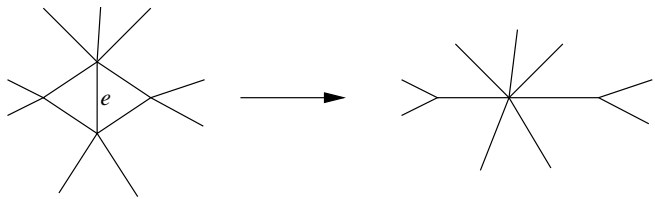}
    \caption{Contraction along an edge}
    \label{fig:contraction}
  \end{center}
\end{figure}

We need some additional notions. A \begriff{contraction} of a regular  
triangulation $T$ along an edge $e$ shrinks $e$ to a vertex and
eliminates the two faces adjacent to $e$, see
Figure~\ref{fig:contraction}. 
The edge $e$ is called \begriff{contractible} if the result
of the contraction is still a regular triangulation. A regular
triangulation $T$ is called \begriff{irreducible} if it does not
contain contractible edges.
The number of vertices of irreducible triangulations is bounded by 
the following result of Nakamoto and Ota~\cite{NakamotoOta95}.
\begin{proposition}\label{thm:finiteirred}
  If $T$ is an irreducible triangulation of a closed surface $F$
  which is not the 
  sphere, then $v(T)\le 270 - 171\chi(F)$. \qed
\end{proposition}

\begin{figure}
  \begin{center}
    \leavevmode
    \epsfig{file=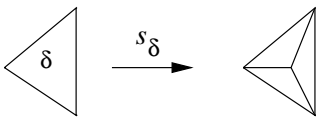}
    \caption{A face subdivision}
    \label{fig:subdivision}
  \end{center}
\end{figure}

Let $\delta$ be a face of a regular triangulation $T$. A \begriff{face 
  subdivision} 
of $T$ along $\delta$ replaces $\delta$ by the cone over its
boundary, see Figure~\ref{fig:subdivision}, 
and the result is denoted
$s_\delta T$. If $\delta$ and $\delta'$ are two faces of $T$, then
$s_\delta T$ and $s_{\delta '} T$ are regularly flip equivalent, which
is easy to see. Let $T_1$ and $T_2$ be regular triangulations and 
$\delta_1$, $\delta_2$ faces of $T_1$, $T_2$. It follows that if $T_1$
and $T_2$ regularly flip equivalent, then so are 
$s_{\delta_1} T_1$ and $s_{\delta_2} 
T_2$. If $T_2$ is obtained from $T_1$ by $m$
successive face subdivisions, then we write $T_2=s^m(T_1)$. The 
notation is ambiguous, but by the preceding remark only up to regular
flip equivalence. After these preliminaries, we can cite a lemma of
Negami~\cite{Negami94}.  
\begin{lemma}\label{lem:1}
  Let $T_1$ and $T_2$ be regular triangulations of $F$. If $T_2$ is
  obtained by contracting some edges of $T_1$, then $T_1$ is
  regularly flip equivalent to  $s^m(T_2)$, with $m= v(T_1) -
  v(T_2)$.\qed   
\end{lemma}

Let $T'$ denote the barycentric subdivision of a singular triangulation
$T$ of  $F$.
\begin{lemma}
  Let $T_1$ and $T_2$ be two singular triangulations of $F$ with
  $v(T_1)=v(T_2)$. Then $T_1''$ and $T_2''$ are regularly flip
  equivalent. 
\end{lemma}
\begin{proof}
It is easy to verify that $T_1''$ and $T_2''$ are regular triangulations
of $F$.
By Proposition~\ref{prop:flips}, we know that $T_1$ and $T_2$ are
related by not necessarily regular flips. 
Let $T_2$ be obtained from $T_1$ by a single  
flip.
Then $T_1'$ can be transformed
into $T_2'$ by the sequence of flips and isotopies that is explicitly
given in Figure~\ref{fig:baryflip}.
\begin{figure}
  \begin{center}
    \leavevmode
    \epsfig{file=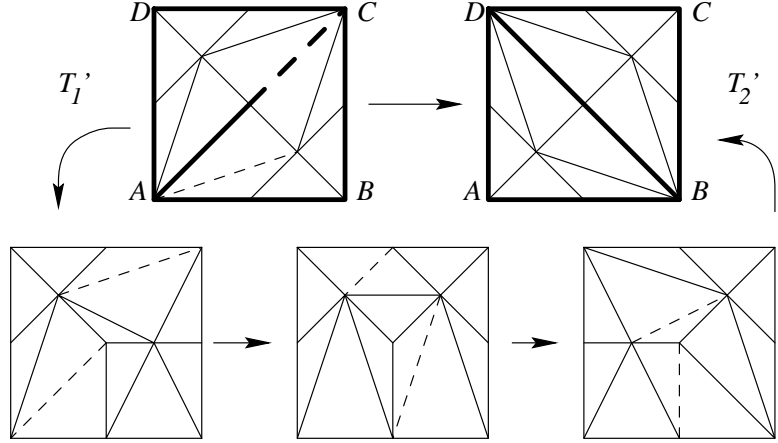}
    \caption{Flips in the barycentric subdivision}
    \label{fig:baryflip}
  \end{center}
\end{figure}
The edges of $T_i$ are drawn bold, and the edges of $T'_i$ under flip  
are dotted. None of these flips introduces a loop.  
It is possible that some flip for $T'_i$ introduces a multiple edge.
This happens only if some of the vertices $A$, $B$, $C$ and $D$
of $T_1$ coincide. We iterate the construction, {i.e.}, we replace each
flip for $T_i'$ by a flip sequence for $T_i''$. 
Since the four vertices of $T'_1$ of each quadrilateral involved in a 
flip  are pairwise distinct, none of these flips introduces a loop or
a multiple edge, thus all flips are regular. \qed
\end{proof}

\begin{corollary}\label{cor:bistellar}
  Let $T_1$ and $T_2$ be two regular triangulations of $F$ with
  $v(T_1)=v(T_2)$. Then $s^m(T_1)$ and $s^m(T_2)$ are regularly flip
  equivalent, with $m=35\left(v(T_1) - 
    \chi(F)\right).$ 
\end{corollary}
\begin{proof}
For any singular triangulation $T$ of $F$, we have $v(T')=v(T)+e(T)+f(T)$. 
Since $2e(T)=3f(T)$, we obtain $f(T)=2(v(T) - \chi(F))$ and  $v(T') =
6v(T)-5\chi(F)$. 
It follows easily  that $v(T_i'')-v(T_i) = m$ for $i=1,2$. 

One obtains $T_i''$ from $T_i$ by $m$ face subdivisions and 
some regular flips, see Figure~\ref{fig:barysubdi} 
\begin{figure}
  \begin{center}
    \leavevmode
    \epsfig{file=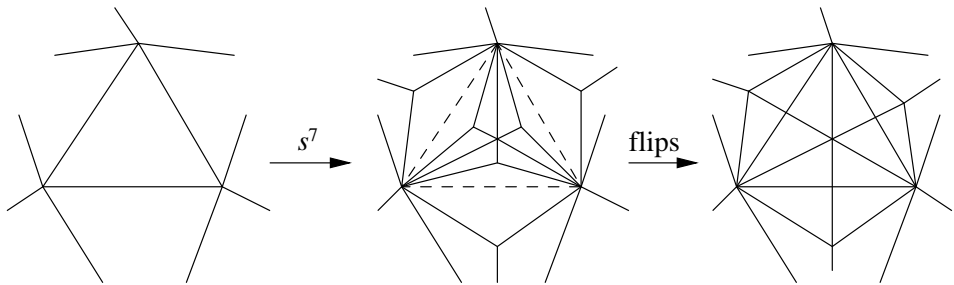,width=\linewidth}
    \caption{How to obtain the barycentric subdivision by flips and
      face subdivisions} 
    \label{fig:barysubdi}
  \end{center}
\end{figure}
for the first
barycentric subdivision. The figure shows the neighbourhood of a face,
and the edges under flip are dotted. So $s^m(T_1)\sim
T_1''\sim T_2''\sim s^m(T_2)$ by the preceding Lemma. \qed
\end{proof}

Now we finish the proof of Theorem~\ref{thm:main}. Let 
$T_1$, $T_2$ be two regular triangulations of $F$ with
$v(T_1)=v(T_2)=M\ge N(F)$ with 
$$N(F)=35\left(\left(270 - 171\chi(F)\right) - \chi(F)\right) = 
                9450 - 6020\chi(F).$$
By contractions along some edges, $T_i$
($i\in\{1,2\}$) can be transformed into an irreducible triangulation
$S_i$. By Lemma~\ref{lem:1}, $T_i$ is regularly flip equivalent to
$s^{M-v(S_i)}S_i$.
By Proposition~\ref{thm:finiteirred} and Corollary~\ref{cor:bistellar},
$s^{N(F)-v(S_1)} S_1$ and $s^{N(F)-v(S_2)} S_2$ are regularly flip
equivalent, and so are also $s^{M-v(S_1)} S_1$ and $s^{M-v(S_2)} S_2$
after further face subdivisions. Therefore also $T_1$ and $T_2$ are
regularly flip equivalent. \QED

\end{document}